\documentclass[12pt,a4paper]{article}
\usepackage{amsmath, amsfonts, amssymb,amsthm}

\begin{document}
\author
{I.Kh. MUSIN}
\title {REPRESENTATION OF INFINITELY DIFFERENTIABLE FUNCTIONS BY DIRICHLET SERIES}
\date{}
\maketitle
\newtheorem{theor}{Theorem}
\newtheorem*{theo}{Theorem A}
\newtheorem{lemma}{Lemma}
\newtheorem{predl}{Proposition}
\newtheorem*{sled}{Corollary}
\newtheorem{sleds}{Corollary}
\newtheorem{note}{Note}
\newcommand{\eps}{\varepsilon}
\begin{abstract}
The problem of representation of elements of
weighted space of infinitely differentiable functions on real line
by exponential series is considered.
\end{abstract}

{\bf 1. Introduction}.
Let  $ \alpha >1 $ and $\psi:{\mathbb R} \to [0,\infty) $
be   a convex function satisfying the conditions:

1.  $ \exists A_{\psi}>0 \
\forall x_1, x_2 \in {\mathbb R} $
$$
\vert \psi (x_1) - \psi(x_2) \vert \le A_{\psi} (1+ \vert x_1 \vert + \vert
x_2 \vert)^{\alpha -1} \vert x_1
- x_2 \vert ;
$$

2. $ \lim \limits_{
x \to \infty}
\displaystyle \frac {\psi (x)} {\vert x \vert} = +\infty $.

Let  ${\cal M}$ be a set of increasing sequences of numbers
$L = (L_n) _{n=0}^{\infty} $ with $L_0 =1$
satisfying the following conditions:

$i_1)$.
$ L_n^2 \le L_{n-1}  L_{n+1}  \
\ \forall n \in {\mathbb N}$;

$i_2)$.
$ \exists \ H_1=H_1(L) > 0, \ H_2 = H_2(L) > 0:
\ \forall \ n \in {\mathbb Z_+} \ \
L_n \ge H_1 H_2^n n!;
$

$i_3)$.  $ \forall \
s > 1 \ \displaystyle \varliminf_{n \to +\infty}
\left( \displaystyle  \frac {L_{[sn]}} {L_n^s} \right)^{\frac 1 n} > 1$;

$i_4)$.
$ \forall \ \delta > 0 \ \exists \
p_{\delta} = p_{\delta}(L)  > 0 \
\exists \ t_{\delta}= t_{\delta}(L) > 1:
\forall \ n \in {\mathbb Z}_+ $
$$
\displaystyle \sup_{m \in {\mathbb N}}
\displaystyle \frac {L_{m+n}}{L_m(1 + \delta )^m}
\le p_{\delta} t^n_{\delta} L_n .
$$

It follows from $ i_4) $ that for any sequence
$ (L_n)_{n=0}^{\infty} $ from ${\cal M}$
\begin{equation}
\lim \limits_{n \to \infty}
{\left
(\displaystyle \frac {L_{n+1}} {L_n} \right)}
^{\frac 1 n} = 1.
\end{equation}

We fix an arbitrary sequence $ M = (M_k)_{k=0}^{\infty} \in {\cal M}$.

Let
$ w(r)= \displaystyle \sup \limits_{k \in {\mathbb Z}_+}
\ln \frac {r^k}{M_k} \ , \  r > 0, \ w(0)=0 $. \
It is easy to see that $w$ is continuous for $ r \ge 0 $.
Also $w(r) = 0 $ for  $ r \in [0, M_1] $. From this and
$i_2)$ it follows that there exists
$A_w >0 $ such that $ w(r) \le A_w r \ , \ r \ge 0. $
Clearly, $w(\vert z \vert)$ is a subharmonic function in the complex plane.

Fix $\sigma >0 $.
Let $ \{{\varepsilon}_m \}_{m=1}^{\infty}$ be  an arbitrary
decreasing to zero sequence of positive  numbers.
Let $ \varphi(x) =
\sup \limits_{y \in {\mathbb R}}
(xy - \psi (y)),  \
{\theta}_m(x) = \exp (\varphi ( x ) - m \ln(1+\vert x \vert)),
\ x \in {\mathbb R}, m \in {\mathbb N} $.
Let
$$
G_m =\{f \in {\cal E}({\mathbb R}): p_m(f) =
\sup_{x \in {\mathbb R}, k  \in {\mathbb Z_+}}
\displaystyle \frac {\vert f^{(k)}(x) \vert}
{(\sigma + {\varepsilon}_m)^k M_k
{\theta}_m(x)} < \infty \}, \
m \in {\mathbb N}.
$$
We let
$G = \bigcap \limits_{m=1}^{\infty}  G_m
$
and endow this vector space with its natural projective limit topology.

Let
$w_m(\vert z \vert) = w ((\sigma +{\varepsilon}_m)^{-1} \vert z \vert),
z \in {\mathbb C}, m \in {\mathbb N}$.
Let
$$
P_m = \left\{f \in H({\mathbb C}):{\Vert f \Vert}_m =
\sup_{z \in  {\mathbb C}}
\displaystyle \frac {\vert f(z) \vert} {\exp(\psi( Im \ z ) +
w_m(\vert z \vert))} < \infty \right\}, m \in {\mathbb N}.
$$
Let $P$ be the union of these normed spaces. The vector space $P$ endowed
with a topology $\tau $ of inductive limit of the spaces  $P_m$ is denoted by
$P_{\tau}$.

For $T \in G^*$ we  define the Fourier-Laplace transform
$ {\hat T }$ of $T$ by
$
{\hat T}(z) = T(e^{-ixz}) , \ z \in {\mathbb C}.
$

Since the sequence $M$ satisfies conditions $i_1), i_2)$ and (1) then the
following theorem holds.
\begin{theo}
The Fourier-Laplace transform establishes topological isomorphism of
the spaces  $G^*$
and
$P_{\tau}$.
\end{theo}

In case
$
\displaystyle \sum \limits_{k =0}^{\infty}
\displaystyle \frac {M_k } {M_{ k +1}}  < \infty
$
theorem A is proved in [2]. In general case theorem A is obtained in [3].

In this paper we announce our result on representation of
functions from $G$ by exponential series.
We pay special attention to examples and properties of the sequences
satisfying conditions $i_1) - i_4)$ and to properties of functions connected
with these sequences.

For weighted spaces of infinitely differentiable functions on real line
similar to $G$ the problem was not considered earlier.

{\bf 2. Examples of sequences belonging to class  ${\cal M}$}.

First introduce the class $ {\cal V} $ of nonnegative convex increasing
functions $ v $ defined on  $ [0, \infty)$ with $ v(0) = 0 $ and such that:

$ {\cal V}1$. $
\exists \ A_v \in {\mathbb R} \ \exists \ B_v \in {\mathbb R}: \forall
x \ge 1 \ \ v(x) \ge x \ln x + A_vx + B_v $;

${\cal V}2$. $  \forall \ s > 1 \
\exists \ {\eta}_s = {\eta}_s(v) > 0 \  \exists \
m_s=m_s(v) \in {\mathbb R}:
\ \forall \ x \ge 0 \ \ v(sx) \ge s v(x) + {\eta}_sx + m_s $;

${\cal V}3$. $  \forall \ \eps > 0
\ \exists \ a_{\eps}=a_{\eps}(v) > 0 \
\exists \ b_{\eps} =b_{\eps}(v) \in {\mathbb R}: \forall \ y \ge 1$
$$
\displaystyle
\sup \limits_{x \ge 1} (v(x+y) - v(x) - \eps x) \le v(y) +
a_{\eps} y + b_{\eps}.
$$

Clearly, for any increasing function $v$ on $[0, \infty] $
satisfying condition ${\cal V}3$  we have
\begin{equation}
\displaystyle \lim \limits_{x \to +\infty} \frac {v(x +1) - v(x)}{x} = 0.
\end{equation}

\begin{predl}
Let $u:[0, \infty) \to [0, \infty) $ be a convex increasing twice continuosly
differentiable function such that:

1. $ \displaystyle \lim \limits_{x \to +\infty}
\frac {u(x +1) - u(x)}{x} = 0;$

2. there exists a constant $ C > 0 $ such that $ u''(x)
\le C x^{-1} $ \ for  $  x \ge 1 $.

Then for any $ \eps \in  (0, C)$ exists a constant
$Q_{\eps}$ such that for all $ y \ge 1 $
$$
\displaystyle \sup \limits_{x \ge 1}
(u(x+y) - u(x) - \eps x)  <  u(y) +
\left(C \ln \displaystyle \frac {2 C} {\eps}  +
\displaystyle \frac {5C}{4}\right) y + Q_{\eps}.
$$

\end{predl}

{\bf Proof}.
Let $ y \ge 1$. Then $ y \in [N, N+1) $ for some  $N \in {\mathbb N}$.
Let $ \eps \in (0, C)$. We shall find the upper estimate of
$ \displaystyle \sup \limits_{x \ge 1} (u(x+y) - u(x) - \eps x - u(y)) $.

>From the first condition on $u$ one can find
a constant $q_{\eps} > 0$  such that
$ u(x+1) < u(x) + \displaystyle \frac {\eps x}{2} + q_{\eps}, x \ge 0 $.
Further, we have
$$
\displaystyle \sup \limits_{x \ge 1}
(u(x+y) - u(x) - \eps x - u(y)) \le
\displaystyle \sup \limits_{x \ge 1}
(u(x+N+1) - u(x) - \eps x - u(N)) \le
$$
$$
\displaystyle \sup \limits_{x \ge 1}
(u(x+N) - u(x) - \displaystyle \frac {\eps x}{2} - u(N+1))
+ \eps N + 2 q_{\eps}.
$$
Note that
$
u(x+N) - u(x) =
\displaystyle \sum \limits_{k= 1}^N
(u(x +k) - u(x +k-1)) \le
\displaystyle \sum \limits_{k= 1}^N u'(x +k);
$
$
u(N+1) =
\displaystyle \sum \limits_{k= 1}^N
(u(k +1) - u(k)) + u(1) >
\displaystyle \sum \limits_{k= 1}^N u'(k).
$
Consequently,
$$
u(x+N) - u(x) - u(N+1) <
\displaystyle \sum \limits_{k= 1}^N (u'(x+k) - u'(k))
=
$$
$$
= \displaystyle \sum \limits_{k= 1}^N
\displaystyle \int \limits_{k}^{x+k} u''(t) \ d t
\le
C \displaystyle \sum \limits_{k= 1}^N
\displaystyle \int \limits_{k}^{x+k}
\displaystyle \frac {d t}{t} =
C \displaystyle \sum \limits_{k= 1}^{N} \ln \left( 1 +
\displaystyle \frac {x}{k}\right).
$$

Thus,
$$
\displaystyle \sup \limits_{x \ge 1}
\left(u(x+N) - u(x) - u(N+1) -
\displaystyle \frac {\eps x}{2} \right)
\le
\displaystyle \sup \limits_{x \ge 1}
\left(C \displaystyle
\sum \limits_{k= 1}^{N}
\ln \left( 1 + \displaystyle \frac {x}{k}\right) -
\displaystyle \frac {\eps x}{2} \right) \le
$$
$$
\sum \limits_{k= 1}^{N}
\displaystyle \sup \limits_{x \ge 1}
\left(C \ln \left( 1 + \displaystyle \frac {x}{k} \right) -
\displaystyle \frac {\eps x}{2N} \right) =
NC \ln \displaystyle \frac {2 C} {e \eps} + C N \ln N -
C \sum \limits_{k= 1}^{N} \ln k +
\displaystyle \frac {\eps (N+1)}{4} \ .
$$

Since
$
\sum \limits_{k= 1}^{N} \ln k \ge
\displaystyle \int \limits_{1}^{N} \ln x \ d x = N \ln N - N + 1
$
then
$$
\displaystyle \sup \limits_{x \ge 1}
\left(u(x+N) - u(x) - u(N+1) -
\displaystyle \frac
{\eps x}{2} \right)
\le
NC \ln \displaystyle \frac {2 C} {\eps}  - C  +
\displaystyle \frac {\eps (N+1)}{4} \ .
$$

Hence, for $ y \ge 1$
$$
\displaystyle \sup \limits_{x \ge 1}
(u(x+y) - u(x) - \eps x - u(y)) <
\left(C \ln \displaystyle \frac {2 C} {\eps}
+ \displaystyle \frac {5C}{4}\right)
y - C + \displaystyle \frac {\eps}{4} + 2 q_{\eps} .
$$
This proves the lemma.

It is easy to see that for arbitrary sequence
$ (L_k)_{k=0}^{\infty} \in {\cal M } $ an increasing function
$v$ on $[0, \infty]$ such that $v(k) = \ln {L}_k, \ k \in {\mathbb Z}_+$,
satisfies conditions ${\cal V}1 - {\cal V}3$. Thus, for arbitrary sequence
$ L = (L_k)_{k=0}^{\infty} \in {\cal M } $
an increasing convex function $v$ on $[0, \infty]$
such that $v(k) = \ln {L}_k, \ k \in {\mathbb Z}_+$,
is in ${\cal V}$. In particular, function $v_L$ such that
$v_L(t k + (1-t) (k+1)) = t \ln {L}_k + (1-t) \ln {L}_{k+1}$,
where $ k \in {\mathbb Z}_+, t \in [0,1]$, is in ${\cal V}$.

Obviously, if increasing function $v$ on $[0, \infty]$
satisfies the condition ${\cal V}3 $ then the sequence
$ (\exp(v(k)))_{k=0}^{\infty} $ satisfies the condition $i_4) $.
Also it is clear that if increasing function $v$ on $[0, \infty]$
satisfies the conditions ${\cal V}2$ and (2) then the sequence
$ (\exp(v(k)))_{k=0}^{\infty} $
satisfies the condition $i_3)$. Thus, for each function
$v \in {\cal V}$ we have $ (\exp(v(k)))_{k=0}^{\infty} \in {\cal M} $.

\begin{predl}
Let $v$ satisfies conditions ${\cal V}1, {\cal V}2$ and conditions of
Proposition 1. Then the sequence $(\exp(v(k)))_{k=0}^{\infty} \in {\cal M}$.
\end{predl}

Now we give some examples of sequences belonging to ${\cal M}$.

1. Consider the function $ v_1(x) = \rho x \ln (x+1), \rho \ge 1, x \ge 0 $.
$v_1 $ is increasing and nonnegative on $[0, \infty), v_1(0) = 0 $.
It is easy to verify that $v_1$ satisfies to conditions ${\cal V}1, {\cal V}2$
and to the first condition of Proposition 1. Since
$
v_1''(x) =
\displaystyle \frac {\rho}{x+1} +
\displaystyle \frac {\rho}{(x+1)^2} > 0$ for $ x \ge 0$ then
$v_1$ is a convex function on $ [0, \infty ) $.
Obviously the second condition of Proposition 1 holds.
By Proposition 2 the sequence
$M^*= ((n+1)^{\rho n})_{n=0}^{\infty}$ is in ${\cal M}$.

Note that for the function
$w^*(r)$ accosiated with the sequence $M^* $
\
$ \rho e^{-1} r^{\frac {1}{\rho}}
- 2 \ln r \le w^*(r) \le \rho e^{-1} r^{\frac {1}{\rho}} $
for $ r > e^{\rho}$

2. Let $ v_2(x) = \rho \ln \Gamma (x + 2),
\rho \ge 1, x \ge 0 $, \ where $\Gamma (x) $ is Euler's Gamma Function.
>From the definition and properties of Gamma Function [5], [6, pp. 755, 763]
it follows that $v_2(0) =0$ and $ v_2 $ is increasing and convex on
$[0, \infty)$ . Using the Stirling's formula it is easy to verify that
the conditions ${\cal V}1, {\cal V}2$ are fulfilled.
The first condition of Proposition 1 is fulfilled obviously. Since
$(\ln \Gamma (x + 2))'' =
\sum \limits_{k= 2}^{\infty}
\displaystyle \frac {1}{(x+k)^2} $
(see, for example, [5], [6, p. 774, formula (28)])
then
$v_2''(x) < \rho
\displaystyle \int \limits_{1}^{\infty}
\displaystyle \frac
{d t} {(x +t)^2} =
\displaystyle \frac {\rho} {x + 1}$. Hence, the second condition
of Proposition 1 is fulfilled too. By Proposition 1 the sequence
$(\Gamma^{\rho} (n+2))_{n=0}^{\infty}$ belongs to ${\cal M}$.

3. For function $ v_3(x) = (x +1) \ln (x + 1) \ {\text arctg} (x +1)$
considered on $[0, \infty) $ we have  $v_3(0) = 0$ and
$v_3'(x) = (\ln (x + 1) + 1) \ {\text arctg} (x +1) +
\displaystyle \frac {(x + 1) \ln (x +1)}{1 + (x+1)^2} > 0$ for $ x \ge 0$.
Hence, $v_3$ is a nonnegative increasing function on $ [0, \infty ) $.
Obviously, condition ${\cal V}_1$ for $v_3$ holds. Since
$
v_3''(x) = \displaystyle \frac {{\text arctg} (x +1)}{x + 1} +
\displaystyle \frac
{2 \ln (x +1)}{(1 + (x+1)^2)^2} +
\displaystyle \frac {2}{1 + (x+1)^2} > 0 $ for  $ x \ge 0$ then
$v_3$ is convex on $ [0, \infty ) $.
Conditions of Proposition 1 are fulfilled since
$
\displaystyle \lim \limits_{x \to +\infty} \frac {v'(x)}{x} = 0
$
and
$v_3''(x) < \displaystyle \frac {6}{x}  $ for $ x \ge 1$.
Next,  for all $ s, x > 1 $
$$
v_3(sx) - s v_3(x) = s x \ln (sx +1) {\text arctg} (sx +1) -
s x \ln (x +1) {\text arctg} (x +1) +
$$
$$
\ln (sx +1) {\text arctg} (sx +1) -
s \ln (x +1) {\text arctg} (x +1)  \ge
$$
$$
\displaystyle \frac {\pi s x}{4} \ln \displaystyle \frac {s+1} {2}  +
\displaystyle \frac
{\pi}{4}
\ln (sx +1) -
\displaystyle \frac
{\pi s}{2} \ln (x +1).
$$
So condition ${\cal V}2$ for $v_3$ is fulfilled. By Proposition 2
the sequence $((n+1)^{(n+1) {\text arctg}(n+1)})_{n=0}^{\infty}$
belongs to ${\cal M}$.

{\bf 3. Auxiliary results}. In this section $v$ is an arbitrary
function in $ {\cal V} $ such that
$ M_k = \exp (v(k)), k \in {\mathbb Z}_+ $. As we know $v$ satisfies
conditions ${\cal V}1 - {\cal V}3$.
Note that conditions  ${\cal V}2$, ${\cal V}3$ impose some conditions
on growth of $v$. For example, from ${\cal V}2$ it follows that for some
$a > 0, b, c  \in {\mathbb R}$ depending on $v$ we have
$ v(x) > a x \ln x + b x + c, \ x \ge 1 $. ${\cal V}3$ implies that for any
$
\eps > 0, x, y \ge 1 $
\begin{equation}
 v(x+y) \le v(x) + \eps x + v(y) +
a_{\eps} y + b_{\eps},
\end{equation}
where the numbers $a_{\eps}>0, b_{\eps}$ depend on $v$ and $\eps$.
In particular, for any $ x \ge 1, \eps > 0$
$ v(2x) \le 2 v(x) + (a_{\eps} + \eps) x + b_{\eps}$. From this inequality
it easily follows that
\begin{equation}
v(x) \le (2v(1) + a_{\eps}+ 2b_{\eps} + \eps)  x +
\displaystyle \frac {(a_{\eps}+ {\eps}) x \ln x}{\ln 2} - b_{\eps}.
\end{equation}

Set
$$
h_v(s) =
\displaystyle \varliminf_{x \to +\infty}
\left(\displaystyle  \frac {v(x)}{x} -
\displaystyle  \frac {v(sx)}{sx} \right), \ s > 0.
$$

\begin{lemma}
Function $ h_v $  has the following properties:

1. for all $ s \in (0, + \infty) \  -\infty < h_v(s) < + \infty $;

2. $h_v(s) > 0 $ for $ s \in (0,1)$ and $h_v(s) < 0 $ for $ s > 1 $;

3. $h_v$ is nonincreasing in $(0,\infty)$;

4. $\displaystyle \lim_{s \to 0, s > 0} h_v(s) = +\infty$;

5. $h_v$ is continuous at the point $s=1$;

6. for any  $ s > 0 \  \ h_v(s) + h_v(s^{-1}) \le 0 $.
\end{lemma}

{\bf Proof}. First note that since $v$ is convex and $ v(0)=0 $
then the function $ \displaystyle  \frac {v(x)}{x} $ is nondecreasing
on $(0, \infty)$. So $h_v(s) \le 0 $ for $s > 1$ and
$h_v(s) \ge 0 $ for $ s < 1 $.

Let $ s > 1 $. Then $ s \in (N, N+1]$ for some  $ N \in {\mathbb N}$ .
>From (3) we have for all $ x \ge (s-N)^{-1}, \eps > 0$
$$
v(sx) \le N v(x) + \eps N x +
\displaystyle  \frac {a_{\eps}Nx(2s -N-1)}{2} + N b_{\eps} +
v((s-N)x).
$$
>From this taking into account that $v(tx) \le tv(x) $ for all
$ t \in [0,1], x \ge 0 $, we get
\begin{equation}
v(sx)  \le sv(x) + \left(\eps  +
\displaystyle  \frac {a_{\eps} (2s - N - 1)}{2}\right)sx + b_{\eps} s,
\end{equation}
for all
$ N \in {\mathbb N}, s \in (N, N+1], x \ge (s-N)^{-1}, \eps > 0$ .
In particular, one can find a constant ${\tilde c_s} > 0 $ such that
for all $ x \ge 0 $
\begin{equation}
v(sx)  \le sv(x) + \left(\eps  +
\displaystyle  \frac {a_{\eps} s}{2}\right)sx + b_{\eps} s +
{\tilde c_s}.
\end{equation}
Using the inequality (6) we obtain
$ h_v(s) \ge -\eps - 0,5a_{\eps} s $
for all $  s > 1, \eps > 0$.

>From the representation
\begin{equation}
h_v(s) = \displaystyle \varliminf_{x \to +\infty}
\left(\displaystyle  \frac {v(s^{-1}x)}{s^{-1}x} -
\displaystyle  \frac {v(x)}{x}\right), \ s > 0,
\end{equation}
and the inequality (6) we have
$ h_v(s) \le \eps + 0,5 a_{\eps}s^{-1} $
for all $ s \in (0,1), \eps > 0$.

Since by the definition $h_v(1) =0$ then the first property is completely
proved.

>From the condition ${\cal V}2$ it follows that $h_v(s) < 0 $ for $ s > 1$.
Using the representation (7) and the condition ${\cal V}2$
we get $h_v(s) > 0 $ for
$ s \in (0,1)$.

The third property of $h_v$ follows from nonincreasing
of $ \displaystyle \frac {v(x)}{x}$ on $ (0, \infty)$.

Since $v$ satisfies the condition ${\cal V}2$ then one can find numbers
$ {\eta}_s(v) > 0, m_s(v)$ such that
$ \forall \ x \ge 0, s> 1 \ \
v(sx) \ge s v(x) + {\eta}_s(v)x + m_s(v) $.
>From this we have for all
$ s > 1, x \ge 0, n \in {\mathbb N} $
$v(s^nx) \ge s^n v(x) + {\eta}_s(v) n s^{n-1} x +
m_s(v) (s^n-1) (s-1)^{-1}$.
Consequently, $ h_v(s^{-n}) \ge s^{-1}{\eta}_sn $ for all $s > 1,
n \in {\mathbb N} $.
>From this and nonincreasing of function $h_v$ we obtain
the fourth property of $h_v$.

We shall prove that function $h_v(s)$ is continuous at $ s=1$.
Let $ \eps > 0$ be arbitrary. Using (5) we have
$ h_v(s) \ge - \eps - a_{\eps}(s-1) $
for $ s \in (1,2) $.
Thus, if $ 0 < s-1 < \min (1, \eps a_{\eps}^{-1})$ then
$ - 2 \eps < h_v(s) - h_v(1) \le 0 $.
Therefore,
$ \displaystyle \lim_{s \to 1, s > 1}h_v(s) = h_v(1) $.
For $ 0,5 < s < 1 $ according to (5)
$v(s^{-1} x) \le s^{-1} v(x) +
\left(\eps  + a_{\eps} (s^{-1} - 1)\right)s^{-1}x +
b_{\eps} s^{-1} $ ,
so $h_v(s) \le  \eps + a_{\eps}(s^{-1} - 1) $.
Therefore, if $ 0 < 1 - s < \min (2^{-1}, \eps (2a_\eps)^{-1})$ then
$0 \le h_v(s) - h_v(1) < 2 \eps $.
Therefore,
$ \displaystyle \lim_{s \to 1, s < 1}h_v(s) = h_v(1) $.
Thus, function $ h_v$ is continuous at point $s= 1$ .

Next, we have
$$
h_v(s^{-1}) =
\displaystyle \varliminf_{x \to +\infty}
\left(\displaystyle  \frac {v(sx)}{sx} -
\displaystyle  \frac {v(x)}{x} \right) =
- \displaystyle {\overline \lim_{x \to +\infty}}
\left(\displaystyle  \frac {v(x)}{x} -
\displaystyle  \frac {v(sx)}{sx}\right) \le
$$
$$
-  \displaystyle \varliminf_{x \to +\infty}
\left(\displaystyle  \frac {v(x)}{x} -
\displaystyle  \frac {v(sx)}{sx} \right) = -h_v(s)
\ ,  s > 0.
$$
>From this we obtain the sixth property of function $h_v$.

Lemma 1 is proved.

\begin{lemma}
The following equality holds
$$
h_v(s) =
\displaystyle \varliminf_{k \to +\infty}
\left(\displaystyle  \frac {v(k)}{k} -
\displaystyle  \frac {v([sk]+1)}{sk}\right)  , \ s >0.
$$
\end{lemma}
{\bf Proof}.
Let $ s >0$. For $ x \in [k, k+1) $,  where $ k \in {\mathbb N}$ ,
we have
$$
\displaystyle  \frac {v(x)}{x} -
\displaystyle  \frac {v(sx)}{sx} \ge
\displaystyle  \frac {v(k)}{k} -
\displaystyle  \frac{v(sk+s)}{sk}=
\displaystyle
\frac {v(k)}{k} - \displaystyle  \frac {v([sk]+1)} {sk}
+
\displaystyle  \frac {v([sk]+1) - v(sk+s)}{sk};
$$
$$
\displaystyle  \frac {v(x)}{x} -
\displaystyle  \frac {v(sx)}{sx} \le
\displaystyle  \frac {v(k+ 1)}{k+1} -
\displaystyle  \frac{v(sk)}{sk} \le
\displaystyle  \frac {v(k)}{k} -
\displaystyle  \frac {v([sk]+1)}{sk}
+
\displaystyle  \frac
{v(sk+1) - v(sk)}{sk} +
$$
$$
+ \displaystyle  \frac {v(k+1) - v(k)}{k+1}
- \displaystyle  \frac {v(k)}{k(k+1)} \ .
$$

Since $v$ satisfies the condition of the form (2) and the inequality (4)
then from the last estimates the assertion of lemma follows.

Thus, according to lemma 2 for any $ v \in {\cal V}$
such that $ \exp (v(k)) = M_k, \ k \in {\mathbb Z}_+ $,
the function $h_v(s)$ coincides with
$$
h(s)=
\displaystyle \varliminf_{k \to +\infty}  (sk)^{-1}
\ln \displaystyle  \frac {M_k^s}{M_{[sk] +1}} \ , \ s > 0.
$$

We set $ l(s) = \exp(h(s)), s > 0$.  From the properties of function
$h$ it follows that function $l$ has the following properties:

1.  for all $ s \in (0, + \infty) \  0 < l(s) < + \infty $;

2. $l(s)$ is continuous at the point $ s=1$;

3. $l(s) > 1 $ for $ s \in (0, 1), \ 0 < l(s) < 1 $ for $ s > 1$;

4. $ \displaystyle \lim_{s \to 0, s > 0} l(s) = +\infty $;

5. $ l(s) l(s^{-1}) \le 1 $.

6. $l$ is a nonincreasing function on  $ (0, \infty)$.

\begin{lemma} For each $ m \in {\mathbb N}$ and $ A > 0 $ there exists a
positive constant $ Q $ such that
$$
w_m(\vert z \vert) + A \ln(1+\vert z \vert) \le
w_{m+1}(\vert z \vert) + Q, \ \ z \in {\mathbb C} .
$$
\end{lemma}

The proof of this lemma is given in [2], [3].

{\bf 4. Weakly sufficient sets for $P$}.
Let $  {\cal K} $ denote a set of all positive continuous functions
$k$ on the complex plane such that for each $ m \in {\mathbb N}$
$$
\displaystyle \sup_{z \in  {\mathbb C}}
\displaystyle \frac
{\exp(\psi( Im \ z ) + w_m(\vert z \vert))} {k(z)}
< \infty .
$$

For each closed subset $S$ of ${\mathbb C}$ that is an uniqueness set
for $P$ we define topologies $ {\tau}_S $ and $ {\mu}_S $
in $P$ in the following manner.
The topology $ {\tau}_S $ is an inductive limit topology of the normed spaces
$$
P_{S,m} = \left\{f \in P:
{\Vert f \Vert}_{S,m} =
\displaystyle
\sup_{z \in  S}
\displaystyle
\frac {\vert f(z) \vert}
{\exp({\psi}(Im \ z) +
w_m(\vert z \vert))} < \infty \right\}, m \in {\mathbb N}.
$$
The topology $ {\mu}_S $ is defined in $P$ with the help of the norms
$$
{\Vert f \Vert}_{S,k} =
\displaystyle
\sup_{z \in  S}
\displaystyle
\frac {\vert f(z) \vert} {k(z)} < \infty , \ k  \in {\cal K}.
$$

Call the subset $S$ sufficient (weakly sufficient)
for $P$ if ${\mu}_{\mathbb C} = {\mu}_S ({\tau}_{\mathbb C} = {\tau}_S)$.

The general arguments [1, chapter 1] show that if $S$ is a sufficient set for
$P$ and  $ \tau ={\mu}_{\mathbb C}$ then every function $ f \in G $ can be
represented as an absolutely convergent integral
$$
f(x) = \int_S
\displaystyle e^{-izx}
\frac
{d \lambda (z)} {k(z)} \ , \ x \in {\mathbb R},
$$
where the complex measure  $\lambda  $ on
$ {\mathbb C}$ is supported by the set $S$ and satisfies the condition
$ \int_{\mathbb C} \vert d \lambda (z) \vert =C_{\lambda} < \infty $,
$ k$ is some function from $ {\cal K}$.
If the sufficient set $S$ is a set of points ${\nu}_j \in {\mathbb C},
\ j = 1, 2, \ldots , $  then
from the integral representation we get the representaion of $f$
in the form of the series
$$
f(x) = \sum_{j=1}^{\infty} c_j e^{-i{\nu}_j x} \,
$$
moreover, from the estimates (see [2], [3]):
$\vert c_j \vert \le
\displaystyle \frac  {C_{\lambda}}{k({\nu}_j)}$
for each
$ j \in {\mathbb N}$, \
$p_m(\exp (-i z x )) \le K_m \exp(\psi ( Im \  z  ) +
w_{m+1}(\vert z \vert ))$
for each $ m \in {\mathbb N},
z \in {\mathbb C}$, where $K_m > 0$ is some constant
independent of $z$,  and lemma 3
it follows that this series absolutely converges in the space $G$.

By the main result of [7] and lemma 3 we have $\tau = {\mu}_{\mathbb C} $.

According to [8], [9] there exists an entire function
${\cal N}(z)$ such that:

1). all the zeros $ \{{\lambda}_j \}_{j=1}^{\infty}$
of  ${\cal N}(z)$ are simple and the discs
$D_j = \{z \in {\mathbb C}: \vert z - {\lambda}_j \vert <
d \} $ are disjoint for some $d > 0$ ;

2). outside the set $ \bigcup \limits_{j=1}^{\infty} D_j$
\begin{equation}
\vert H_D(z) + w(\sigma^{-1} \vert z \vert)
- \ln \vert {\cal N}(z) \vert \vert
\le A \ln(1 + \vert z \vert ) + C_0 ,
\end{equation}
where $A, C_0$ are some positive numbers.

\begin{theor}
The set ${\tilde S}= \{{\lambda}_j \}_{j=1}^{\infty}$ of zeros of
${\cal N}$ is a weakly sufficient set for $P_{\tau}$.\end{theor}

\begin{theor}
Every function $f \in G $ can be represented in the form of a series
$$
f(x) = \sum_{j=1}^{\infty} c_j e^{-i{\lambda_j} x},
$$
absolutely convergent in $G$.
\end{theor}

The proof of Theorem 1 is based on the representation of entire functions
in the space $P$ by Lagrange series, on the key
\begin{lemma}
For all $s > 0, \delta \in (0,1)$ there exists a constant
$Q(s, \delta)\ge 0$ such that
 $ s w(r)
\le w \left(\displaystyle \frac  r {l(s)(1 - \delta)}\right) +
Q(s,\delta) $  for all $ r \ge 0 $
\end{lemma}
and will be given in [4].
\pagebreak

\begin{center}
Bibliography
\end{center}

1. L. Ehrenpreis.  {\it Fourier analysis in several complex variables}.

Wiley -- Interscience, New York, 1970.

2. I. Kh. Musin. Paley-Wiener type theorem for a weighted space of

infinitely differentiable functions. Izv. Akad. Nauk SSSR. Ser. Mat.,

64:6 (2000),    181-204.

3. I. Kh. Musin. On  the Fourier-Laplace  transform of functionals
on a

weighted space of infinitely differentiable functions. Pbb:

funct-an@xxx.lanl.gov  N 9911067.

4. I. Kh. Musin.  On representation of infinitely differentiable

functions by Dirichlet seies. (submitted to "Matematicheskie zametki")

5. F. W. J. Olver. Introduction to asymptotics and special functions.

Academic Press, NY and London, 1974.

6. G. M. Phihtengoltc. Course of differential and integral calculus.

V. II. M.- L.: GITTL, 1948.

7. V. V. Napalkov. On comparison of the topologies
in some spaces of

entire functions.  Dokl. Akad. Nauk SSSR, 264:4 (1982), 827-830.

8. R. S. Youlmukhametov. Approximation of subharmonic functions.

Math. Sb., 1984. 124:3(1984), 393-415.

9. R. S. Youlmukhametov. Approximation of subharmonic functions.

Analysis Mathematica, 11:3 (1985), 257-282.

\end{document}